\documentclass{article}
\usepackage{amsfonts,amsmath}
\title{$\mathbb{Z}_{2}^{2}$-cordiality of $K_{n}$ and $K_{m,n}$}
\author{Adrian Riskin
\\
Department of Mathematics
\\
Mary Baldwin College
\\
Staunton, Virgina 24401
\\
ariskin@mbc.edu}

\begin{document}

\maketitle

\abstract{We prove that $K_{n}$ is $\mathbb{Z}_{2}^{2}$-cordial if and only if $1 \leq n \leq 3$
and that $K_{m,n}$ if and only if it's false that $m=n=2$}
\section{Introduction}
Cahit [1] introduced the notion of cordial graphs.  For an introduction and comprehensive
survey see [2]. This notion was generalized by Hovey [3] to $A$-cordiality.  A graph $G$ is 
$A$-cordial where $A$ is an additive abelian group and $f:V(G) \rightarrow A$ is a labeling
of the vertices of $G$ with elements of $A$ such that when the edges of $G$ are labeled by 
the induced labeling $f_{e}:E(G) \rightarrow A$ by $f(ab)=f(a)+f(b)$ then for each
$a,b \in A$ we have $|f^{-1}(a)-f^{-1}(b)| \leq 1$ and likewise for $f^{-1}_{e}$.  Hovey calls
a graph $k$-cordial if $A=\mathbb{Z}_{k}$.  Note therefore that $2$-cordiality is ordinary 
cordiality.  Heretofore research on $A$-cordiality has been limited to research on $k$-cordiality.
Herein we study the $\mathbb{Z}_{2}^{2}$-cordiality of $K_{n}$ and $K_{m,n}$.  In addition 
we have a few minor results and questions on the $\mathbb{Z}_{2}^{2}$-cordiality of trees, since
 trees have been such a focus of interest in research on $k$-cordiality.  We call a vertex labeling 
 {\it balanced} if it meets the above criterion for vertices of $A$-cordial labelings, and likewise with 
 an edge labeling. 
 
 \section{ The $\mathbb{Z}_{2}^{2}$-cordiality of $K_{n}$}
 
 Cahit [1] has shown that $K_{n}$ is cordial if and only if $1 \leq n \leq 3$.  We have 
 the following analogous result:
 
 \bigskip
 
 \noindent \textbf{Theorem 1}: $K_{n}$ is $\mathbb{Z}_{2}^{2}$-cordial if and only if $1 \leq n \leq 3$.
 
 \bigskip
 
 \noindent \textbf{Proof}: For $1 \leq n \leq 3$ $K_{n}$ is easily seen to be $\mathbb{Z}_{2}^{2}$-cordial.
 Now, suppose that $n \geq 4$ and (first of all) that $n \equiv 0 \pmod 4$.  If $K_{n}$ were to 
 be cordial we'd have $\frac{n}{4}$ vertices labeled with each element of $\mathbb{Z}_{2}^{2}$.  This yields
 $$4 {\frac{n}{4} \choose 2}$$
 edges labeled 00 and
 $$2\left(\frac{n}{4}\right)^{2}$$
 edges labeled 11.  If $K_{n}$ is to be $\mathbb{Z}_{2}^{2}$-cordial we must have 
 $$1 \geq \left | 2 \left (\frac{n}{4}\right)^{2}-4{\frac{n}{4}}\right| = \left |\frac{n}{2}\right|$$
 so that $2 \geq n$, an impossibility.
 
 If $n\geq 4$ and $n \equiv 1 \pmod 4$ we have four distinct vertex labelings of $K_{n}$ compatible with 
 $\mathbb{Z}_{2}^{2}$-cordiality determined by which element of $\mathbb{Z}_{2}^{2}$ labels $\frac{n+3}{4}$
 vertices.  In each of these four cases we have
 $$3{\frac{n-1}{4} \choose 2} + {\frac{n+3}{4} \choose 2}$$
 edges labeled 00 and 
 $$\left(\frac{n+3}{4}\right)\left(\frac{n-1}{4}\right) + \left (\frac{n-1}{4}\right)^{2}$$
 edges labeled 11.  In order for $K_{n}$ to be $\mathbb{Z}_{2}^{2}$-cordial we must have
 $$1 \geq \left | 3{\frac{n-1}{4} \choose 2} + {\frac{n+3}{4} \choose 2} - 
 \left(\frac{n+3}{4}\right)\left(\frac{n-1}{4}\right) + \left (\frac{n-1}{4}\right)^{2} \right |
 = \frac{1}{8} |5n-4|$$
 Hence
 $$8 \geq 5n-4 \geq -8$$
 or
 $$\frac{12}{5} \geq n \geq -\frac{4}{5}$$
 which is again impossible.  The other two cases are handled similarly. \hfill $\square$
 
 \bigskip
 
 We say that an $A$-cordial labeling of graph $G$ is {\it A-monocordial} if each element appears either once or 
 not at all upon each vertex and each edge of $G$.  Clearly if there is an $A$-monocordial labeling of $G$ and
 a monomorphism $\phi:A \rightarrow B$ then there is an $B$-monocordial labeling of $G$ as well.  Since 
 there is a $\mathbb{Z}_{2}^{2^{n-1}}$-monocordial labeling of $K_{n}$ obtained by numbering
 the vertices of $K_{n}$ from 0 through $n-1$ and labeling vertex $i$ with the element of 
 $\mathbb{Z}_{2}^{2^{n-1}}$ with a 1 in the $i^{th}$ place and 0's elsewhere, and labeling vertex 0 with the 
 identity element of $\mathbb{Z}_{2}^{2^{n-1}}$.  Thus for each $n$ there is a value of $m$ such that 
 $K_{n}$ is $\mathbb{Z}_{2}^{p}$-cordially labeled for $p \geq m$.  This raises the question of finding 
 the minimum such $m$ for each $n$.  So far, then, it is known that for $n=4$ the minimum value of $m$ is 3, but 
 no other values are known.
 
 \section{$\mathbb{Z}_{2}^{2}$-cordiality of trees}
 
 Cahit [1] showed that all trees are 2-cordial.  Hovey [3] showed that all trees are $k$-cordial for 
 $k=3, 4, 5$ and conjectures that all trees are $k$-cordial for all $k$.  On the other hand it is easily 
 seen that the path with 5 vertices is not $\mathbb{Z}_{2}^{2}$-cordial by an examination of cases.  This 
 suggests the seemingly difficult problem of determining exactly which trees are $\mathbb{Z}_{2}^{2}$-
 cordial and which are not.
 
 \bigskip
 
 \noindent \textbf{Theorem 2}: The $n$-star $K_{1,n}$ is $\mathbb{Z}_{2}^{2}$-cordial.
 
 \bigskip
 
 \noindent \textbf{Proof:} Label the $n$ end vertices with elements of $\mathbb{Z}_{2}^{2}$ so that the labeling
 is balanced, with label 00 in one of the smaller sets of labels if any such there are.  Label the central vertex
 00.  Clearly the induced edge labeling is balanced as well. \hfill $\square$
 
 \bigskip
 
 This theorem also raises the question of which complete bipartite graphs are $\mathbb{Z}_{2}^{2}$-cordial.
 This question is answered in the next section.
 
 \section {$\mathbb{Z}_{2}^{2}$-cordiality of $K_{m,n}$}
 
 Our main result here is that $K_{m,n}$ is $\mathbb{Z}_{2}^{2}$-cordial except in the case of $K_{2,2}$.
 
 \bigskip
 
 \noindent \textbf{Theorem 3:} If one of $m$ and $n$ is a multiple of 4 then $K_{m,n}$ is 
 $\mathbb{Z}_{2}^{2}$-cordial.
 
 \bigskip
 
 \noindent \textbf{Proof:} Suppose without loss of generality that 4 divides $m$.  Let $A$ and $B$ be the 
 $m-$ and $n-$partite sets of $K_{m,n}$ respectively.  Label the vertices of $A$ with $\frac{m}{4}$ copies
 of the elements of $\mathbb{Z}_{2}^{2}$.  Label the vertices of $B$ from $\mathbb{Z}_{2}^{2}$ in a 
 balanced manner.  Then for $b \in B$ the labels of the edges incident with $b$ comprise $\frac{m}{4}$
 copies of each element of $\mathbb{Z}_{2}^{2}$ since the edges incident with $b$ and with each of the 
 $\frac{m}{4}$ copies of $\mathbb{Z}_{2}^{2}$ have labels which constitute a coset of $\mathbb{Z}_{2}^{2}$ in
 itself, and are therefore distinct.  Hence the edges are labeled with $\frac{mn}{4}$ copies of each of the 
 elements of $\mathbb{Z}_{2}^{2}$, which produces a $\mathbb{Z}_{2}^{2}$-cordial labeling. \hfill $\square$
 
 \bigskip
 
 \noindent \textbf{Lemma 1:} $K_{2,2}$ is not $\mathbb{Z}_{2}^{2}$-cordial.
 
 \bigskip
 
 \noindent \textbf{Proof:} Examine cases. \hfill $\square$
 
 \bigskip
 
 \noindent \textbf{Theorem 4:} $K_{2,n}$ is $\mathbb{Z}_{2}^{2}$-cordial for $n \geq 3$.
 
 \bigskip
 
 \noindent \textbf{Proof:} It is easily seen that $K_{2,3}$ is $\mathbb{Z}_{2}^{2}$-cordial by construction.
 Also, $K_{2,4}$ is cordial by Theorem 3.  Hence we may assume $n \geq 5$.  Suppose $n \equiv i \pmod 4$ with
 $0 \leq i \leq 3$.  Label the $n$-set of vertices with $\frac{n-i}{4}$ copies of $\mathbb{Z}_{2}^{2}$. Label
 the vertices of the 2-set 11 and 10.  If $i=0$ we're done.  If $i=1$ label the remaining vertex of the 
 $n$-set 00.  If $i=2$ label the two remaining vertices of the $n$-set 00 and 10.  If $i=3$ label the three 
 remaining vertices of the $n$-set 00, 10, and 01.  In each case it is easily checked that the result is
 a $\mathbb{Z}_{2}^{2}$-cordial labeling. \hfill $\square$
 
 \bigskip
 
 \noindent \textbf{Theorem 5}: $K_{3,n}$ is $\mathbb{Z}_{2}^{2}$-cordial for $n \geq 3$.
 
 \bigskip
 
 \noindent \textbf{Proof:} $K_{3,3}$ is easily seen to be $\mathbb{Z}_{2}^{2}$-cordial by construction.  If
 $n \geq 4$, suppose $n \equiv i \pmod 4$ with $0 \leq i \leq 3$.  If $i=0$ we're done by Theorem 3.  If $i=1$
 label the last vertex of the $n$-set with 10 and the three vertices of the 3-set with 00, 01, and 11.  Then 
 we have $\frac{n-1}{4}$ edges labeled 00 and $\frac{n+3}{4}$ labeled with each of 10, 01, and 11.  Thus the 
 labeling is $\mathbb{Z}_{2}^{2}$-cordial.  If $i=2$ we label the last two vertices of the $n$-set with 10 and
 00, and if $i=3$ we label the last three vertices 10, 00, and 11.  Count in the same way as above to verify 
 that these indeed provide $\mathbb{Z}_{2}^{2}$-cordial labelings. \hfill $\square$
 
 \bigskip
 
 \noindent \textbf{Theorem 6}: $K_{m,n}$ is $\mathbb{Z}_{2}^{2}$-cordial except for the case where $m=n=2$.
 
 \bigskip
 
 \noindent \textbf{Proof:} By above work we may assume that $m \geq n \geq 5$.  Suppose $m \equiv i \pmod 4$ with
 $0 \leq i \leq 3$ and $n \equiv j \pmod 4$ with $0 \leq j \leq 3$.  If either $i$ or $j$ is 0 then we're done.
 Label $m-i$ vertices of the $m$-set with $\frac{m-1}{4}$ copies of  $\mathbb{Z}_{2}^{2}$ and $n-j$ of the 
 vertices of the $n-set$ with $\frac{n-j}{4}$ copies of $\mathbb{Z}_{2}^{2}$.  Note that among the 
 $(m-i)(n-j)$ labeled edges in this part of the graph exactly one fourth are labeled with each element of 
 $\mathbb{Z}_{2}^{2}$ and that the same is true of the vertices.  The unlabeled portion of the graph is 
 isomorphic to $K_{i,j}$.  Unless $i=j=2$, label this $K_{i,j}$ $\mathbb{Z}_{2}^{2}$-cordially, as we know is 
 possible from above work.  Then one fourth of the $j(m-i)+i(n-j)$ edges running from this $K_{i,j}$ isomorph to 
 the $K_{m-i,n-j}$ isomorph are labeled with each of the elements of $\mathbb{Z}_{2}^{2}$, and the edge labelings
 within the $K_{i,j}$ are balanced, making the whole thing a $\mathbb{Z}_{2}^{2}$-cordial labeling.  If
 $i=j=2$ label the vertices of the $K_{2,2}$ isomorph with the four elements of $\mathbb{Z}_{2}^{2}$.  The same sort 
 of count shows that the entire graph ends up $\mathbb{Z}_{2}^{2}$-cordially labeled. \hfill $\square$
 
 \section*{References}
 
 \begin{enumerate}
 \item Cahit, I.; Cordial graphs: a weaker version of graceful and harmonious
graphs.  Ars Combin. 23(1987) 201-207.

\item Gallian, J. A.; A dynamic survey of graph labeling.  Electronic J.
Combin. DS6. http://www.combinatorics.org/Surveys/index.html

\item M. Hovey, A-cordial graphs, Discrete Math., 93 (1991) 183-194.
 
 \end{enumerate}

\end{document}